# Higher Derivatives of the Falling Factorial and Related Generalizations of the Stirling and Harmonic Numbers

By Steven S. Poon

January 12, 2014

**Introduction**

This work started as a simple exercise to calculate solutions for the higher derivatives of the falling factorial function evaluated at integral values:

$$F_n^{(l)}(m) = \frac{d^l}{dx^l}\{x^{\underline{n}}\}\bigg|_{x=m} \quad \begin{array}{l} n \in \mathbb{N} \\ m \in \mathbb{Z} \\ x \in \mathbb{C} \end{array} \quad (1)$$

There is no doubt that this derivative can be calculated easily with elementary methods for small fixed values of $l$. As a rule, deriving a general solution for a variable value of $l$ is a more difficult task.

The general solution was first published by Koutras [11], who gave the name "non-central Stirling numbers" to the following quantity.

$$s_m(n, l) = \frac{1}{l!} F_n^{(l)}(m) \quad (2)$$

His method was to first find out, through sound reasoning, the recursive relationship of these numbers. He was then able to manipulate the recursive relationship to derive a solution of these numbers as a weighted sum of the Stirling numbers of the first kind. It is straightforward to apply these formulae to calculate numerical solutions for (1).

Another development came from Broder [12]. His work started by examining the combinatorial properties of a set of quantities that he named the *r*-Stirling numbers. He defined the *r*-Stirling numbers of the first kind as "the number of permutations of the set {1, …, n} having *m* cycles, such that the numbers 1, 2, …, *r* are in distinct cycles". He was able to derive a large number of interesting identities involving these numbers. Among these identities is a solution for $F_n^{(l)}(m)$ using *r*-Stirling numbers, which happens to be the higher derivatives of the generating function:

$$(x+r)^{\overline{n}} = \sum_{k} \begin{bmatrix} n+r \\ k+r \end{bmatrix}_r x^k \qquad (3)$$

Note that the square bracket with the subscript is Broder's notation for the *r*-Stirling numbers.

The focus of the current work is to describe a brute-force solution for (1) through repeated applications of the product rule (also known as Leibniz rule) of differentiation. Because of the large number of terms that result from such an approach, methods for annotating the results will be described. This article will also attempt to fill in a few insights that might not have been emphasized by previous authors who approached this question through different angles. It will be shown that many different quantities, including certain generalized Stirling numbers, harmonic numbers, and symmetric polynomials, are, in fact, a coherent whole.

## Some Historical Notes

In the early 20[th] century, Professor Niels Nielsen was credited as the first person who coined the term "Stirling numbers" [4-6], in recognition of James Stirling's pioneering work on these quantities, which was published in his treatise in 1730 [3]. Stirling appeared to be the first to make use of one of the most recognizable properties of these numbers:

$$x^{\underline{n}} = \sum_{k \geq 0} s(n,k) x^k \qquad \begin{matrix} n \in \mathbb{N} \\ x \in \mathbb{C} \end{matrix} \qquad (4)$$

$$x^n = \sum_{k \geq 0} S(n,k) x^{\underline{k}} \quad \begin{matrix} n \in \mathbb{N} \\ x \in \mathbb{C} \end{matrix} \tag{5}$$

Stirling himself did not create a dedicated symbol or notation for the Stirling numbers, and simply used the letters *A*, *B*, *C*, etc. to denote different Stirling numbers of the same row. Since then, the usefulness of these numbers has been discovered in often unexpected ways by different groups of researchers, resulting in a lack of standardization of the notation. Knuth had a good discussion of this topic in his work "Two Notes on Notations" [1]. However, this work will be using the notation for Stirling numbers used by John Riordan is his book *Combinatorial Identities* [2]. Stirling numbers of the first and second kinds will be annotated as $s(n, k)$ and $S(n, k)$, respectively. This is because Knuth's notation is associated with a "sign-less" definition of the Stirling numbers of the first kind, which is not compatible with the results to be shown later. However, Knuth's notations for the falling factorials and Iverson's convention will be used in this paper.

It should be noted that equation (1) is a clear statement that $s(n, k)$ are the coefficients of the falling factorial when expanded in standard polynomial form. As a result, these Stirling numbers are destined to be intimately related to the concept of symmetric functions.

Viète [7] and Girard [8] are often credited as the pioneers of symmetric polynomials. Even though they did not use the exact term, nor did they note the symmetric properties of these functions, they did describe, without proof, the form of the elementary symmetric polynomials. In modern notation, these functions can be written as:

$$e_k(x_1, x_2, \cdots, x_n) = \begin{cases} 1 & k = 0 \\ \sum_{1 \leq t_1 < t_2 < \cdots < t_k \leq n} \left( \prod_{j=1}^{k} x_{t_j} \right) & 1 \leq k \leq n \end{cases} \tag{6}$$

Furthermore, it was Girard who stated clearly that functions of this form are the coefficients of a polynomial with roots $x_1, x_2, \ldots, x_n$ (this relationship is often called Viète's formulas). Since the falling factorial is a polynomial with roots 0, 1, …, (*n*

– 1), one may evaluate the elementary symmetric polynomial at these values to give:

$$|s(n, n-k)| = e_k(0,1,\cdots,n-1) \quad n \geq k \geq 1 \tag{7}$$

It is thus clear that one way to generalize the Stirling numbers is to manipulate the limits in equation (7). It will be shown below (as Broder has shown previously [12]) that such generalized quantities are natural ways to annotate the solutions for (1) when the variable $x$ is an integer.

### Differentiating the Falling Factorial Using the Product Rule

The first derivative of the falling factorial is known to be studied by Isaac Newton [13-14]. The exercise itself is an elementary application of the Leibniz rule of differentiation for any finite and fixed value of $n$ (although Newton himself might have used other methods). When $n$ is considered a variable value, a rigorous proof may be achieved using an inductive argument. This proof is quite simple and will not be repeated here.

A more general result applies to the derivative of a product of $n$ factors, each in the form $(x - k_j)$, where $k_j$ is not equal to $k_{j'}$ if $j$ is not equal to $j'$. Differentiating such a quantity results in $n$ terms, with each term being the same as the original product with one factor removed. There are $n$ ways to remove one factor from $n$ factors, which correspond nicely with the $n$ terms generated by the differentiation.

$$\frac{d}{dx}\left\{\prod_{j=0}^{n-1}(x-k_j)\right\} = \sum_{r=0}^{n-1}\left[\frac{1}{x-k_r}\prod_{j=0}^{n-1}(x-k_j)\right] \tag{8}$$

Applying this result to the falling factorial function gives:

$$F_n'(x) = \begin{cases} 0 & n = 0 \\ x^{\underline{n}} \sum_{k=0}^{n-1} \dfrac{1}{x-k} & n \geq 1 \end{cases} \quad x \in \mathbb{C} \wedge x \notin \mathcal{N}_n \tag{9}$$

This leads to well-known relationships with Euler's digamma function and harmonic numbers.

$$F_n'(x) = [\varphi(x+1) - \varphi(x-n+1)]x^{\underline{n}} \tag{10}$$

$$F_n'(m) = (H_m - H_{m-n})m^{\underline{n}} \quad m \in \mathbb{Z} \wedge m \geq n+1 \tag{11}$$

The standard treatment of the problem of calculating the derivatives of the falling factorial function tends to stop at the first derivative. It is certainly a simple exercise to extend the calculation to any fixed small values of $n$. However, one can already see the beginning of several interesting difficulties, which will be described below.

## The Divide-by-Zero Problem for Calculating the Derivative of the Falling Factorial Function

The equation (9) is in an unfortunate form. The left hand side, which is the derivative of the falling factorial power, is defined for all complex values. The right hand side, on the other hand, is not well-defined when $x$ is an integer ranging between 0 and $(n-1)$. Even if one argues that the equation is still correct when taken as a limit, the exact form of the equation still cannot be used directly for numerical calculation purposes due to the lack of clear directions on how to cancel out the pole.

This situation appears to be caused by the fact that the standard treatment has a slightly incorrect interpretation of the nature of the derivative of the falling factorial. Fundamentally, differentiating the falling factorial of order $n$ creates $n$ terms of falling factorials with one missing factor for each term. It is true that the falling factorial with one missing factor may be approximated by dividing out the factor to be removed (this is the strategy taken in equation (9)). However, this approach runs into the divide-by-zero problem. A more correct approach is to remove the unwanted factor, instead of dividing it out.

$$F'_n(x) = \begin{cases} 0 & n = 0 \\ \sum_{k=0}^{n-1} \left[ \prod_{j \in \{i : i \in \mathbb{N}_n \wedge i \neq k\}} (x - j) \right] & n \geq 1 \end{cases} \quad x \in \mathbb{C} \qquad (12)$$

The simplicity of this fix appears to undermine the assertion that the divide-by-zero problem is a significant issue. One immediate negative effect, however, is that the connection to harmonic numbers has become less apparent with this notation. Also, as one will see, the number of missing factors will also increase as the order of the differentiation operator increases. The notation used in equation (12) will become increasingly unwieldy. The exact algorithmic interpretation of the equation will also become less clear. It will be desirable to create an alternative way to annotate a falling factorial with missing factors in a way that the steps for numerical calculation are unambiguous.

## Introducing the Falling Factorial with Missing Factors

For the purpose of this article, the falling factorial of order *n* and *l* missing factors is defined as follows.

$$\Theta(n, \vec{k}_l; x') = \lim_{x \to x'} \left\{ x^{\underline{n}} \prod_{j=1}^{l} \frac{1}{x - k_j} \right\} \qquad (13)$$

$$\vec{k}_l = \langle k_j | j \in \mathbb{Z} \wedge 1 \leq j \leq l \rangle \quad k_j \in \mathbb{Z} \wedge 0 \leq k_j \leq n - 1 \qquad (14)$$

The *l*-tuple <$k_1$, ..., $k_l$> contains the missing factors that are to be removed from the falling factorial. All members of the *l*-tuple must be integers ranging from 0 to (*n* − 1). Strictly speaking, these members do not have to be unique. However, if the members of the *l*-tuple do not repeat, equation (13) is a simple polynomial of order (*n* − *l*). If there are repeating factors, equation (13) is a polynomial added to a rational function. Either way, the polynomial portion, including its coefficients, may be calculated through long division.

For a concrete example, consider the case *l* = 1. Performing the long division on equation (9) yields the following result.

$$\Theta(n, k_1; x) = \sum_{j=0}^{n-1} \vartheta(n, k_1, j) x^j \tag{15}$$

$$\vartheta(n, k_1, j) = \begin{cases} s(n, j+1) & (k_1 = 0) \\ 0 & (k_1 > 0, j = 0) \\ -\sum_{i=1}^{j} s(n, i) k_1^{i-j-1} & (k_1 > 0, j > 0) \end{cases} \tag{16}$$

One may easily verify this through direct evaluation.

These coefficients retain the defining recursive relationship of the Stirling numbers of the first kind, as shown in the equation below. This relationship may be confirmed easily by plugging it back into equation (16).

$$\vartheta(n+1, k_1, j) = \vartheta(n, k_1, j-1) - n\vartheta(n, k_1, j) \tag{17}$$

As a result, these coefficients may be considered generalizations of the Stirling numbers of the first kind. Because the same recursive relationship is followed, the only difference must be the initial values.

In addition to equation (17), the coefficients $\vartheta(n, k_1, j)$ also observe the following identities. Both of these identities are easily proven through direct evaluation of equation (16).

$$\vartheta(n, k_1, n-1) = 1 \tag{18}$$

$$\vartheta(n, k_1, j) - k_1 \vartheta(n, k_1, j+1) = \vartheta(n, 0, j) = s(n, j+1) \tag{19}$$

One may continue to calculate the coefficients $\vartheta(n, <k_1, \ldots, k_l>, j)$ for higher and higher fixed values of $l$ through the long division method. For example, the coefficients for $\vartheta(n, <k_1, k_2>, j)$ may be obtained by dividing equation (15) by $(x - k_2)$. It is harder to obtain a general solution for a variable value of $l$. The obvious way to obtain such a solution is to calculate the solutions for the first few fixed values of $l$, form a conjecture, and prove the conjecture through induction. The following is such a conjecture.

$$\vartheta(n, \vec{k}_l, j) = \sum_{i=0}^{n-j-l} \vartheta(n, \vec{k}_{l-1}, i+j+1) k_l^i \qquad (20)$$

If one assumes the above is correct for a particular value of $l$, it is then easily shown to be true for $(l + 1)$ by going through the long division exercise.

This quantity follows similar identities that $\vartheta(n, k_1, j)$ also follows. In particular, the following identities are all true.

$$\vartheta(n, \vec{k}_l, n - s) = 1 \qquad (21)$$

$$\vartheta(n, \vec{k}_l, j) - k_l \vartheta(n, \vec{k}_l, j+1) = \vartheta(n, \vec{k}_{l-1}, j+1) \qquad (22)$$

$$\vartheta(n+1, \vec{k}_l, j) = \vartheta(n, \vec{k}_l, j-1) - n\vartheta(n, \vec{k}_l, j) \qquad (23)$$

$$\vartheta(n+1, \vec{k}_l, j+1) + (n - k_l)\vartheta(n, \vec{k}_l, j+1) = \vartheta(n, \vec{k}_{l-1}, j+1) \qquad (24)$$

Equation (21) should be evident since the falling factorial with missing factors is made up of factors in the form of $(x - k_j)$. The coefficient of the highest-order term is necessarily equal to 1. Equation (22) is easily obtained from equation (21) by changing the indexing scheme of the summation such that $i$ is shifted by one.

$$k_l \vartheta(n, \vec{k}_l, j+1) = \sum_{i=0}^{n-j-l} \vartheta(n, \vec{k}_{l-1}, i+j+1) k_l^i - \vartheta(n, \vec{k}_{l-1}, i+1) \qquad (25)$$

This then immediately reduces to the same form as equation (22).

Equation (23) may also be proven through induction. Suppose one assumes that it is true for a certain value of $l$. Expanding equation (21) using the assumed recursive relationship gives:

$$\vartheta(n+1, \vec{k}_l, j) = \sum_{i=0}^{n-j-l+1} \vartheta(n, \vec{k}_{l-1}, i+j) k_l^i \\ - n \sum_{i=0}^{n-j-l+1} \vartheta(n, \vec{k}_{l-1}, i+j+1) k_l^i \tag{26}$$

This reduces to the desired identity (23) because the last term in the second summation is zero.

Finally, equation (24) is just a combination of (22) and (23).

Equation (22) is particularly interesting because it shows that the coefficients ϑ(n, <$k_1$, …, $k_l$>, j) follow the same recursive relationship as that of the Stirling numbers of the first kind, no matter how high the value of $l$ is. Thus these coefficients may be seen as a generalization of the Stirling numbers of the first kind. This is the kind of generalization where the recursive relationship remains unchanged, while the initial values are changed to a different set.

Equation (22) is also remarkable because it allows one to pick any "row" of values ϑ(n, <$k_1$, …, $k_l$>, j) with n and <$k_1$, …, $k_l$> fixed, and proceed to using those values to create a complete matrix of values for ϑ(n, <$k_1$, …, $k_{l+1}$>, j) with all possible values of $k_{l+1}$. As a result, one may calculate the set of coefficients for any falling factorial with missing factors through a purely recursive algorithm, starting from the Stirling numbers of the first kind.

With the concept of the falling factorial with missing factors now fully understood, one may rewrite the solution of the first derivative of the falling factorial function as:

$$F_n'(x) = \begin{cases} 0 & n = 0 \\ \sum_{k=0}^{n-1} \Theta(n, k; x) & n \geq 1 \end{cases} \quad x \in \mathbb{C} \tag{27}$$

### Enumerating All Possible Sets of Missing Factors

One more piece of background is needed before plunging into the discussion of solving the higher-order derivatives of the falling factorial function. This is related to the need of coming up with a good way to annotate all the possible terms that are created when the falling factorial power is processed by the differentiation operator repeatedly.

Consider the power set, or the set of all subsets, of $\mathcal{N}_n$. It may be used to define the following set:

$$C^*_{n,l} = \{\mathcal{X} : \mathcal{X} \in \mathcal{P}(\mathcal{N}_n) \wedge |\mathcal{X}| = l\} \quad 0 \leq l \leq n \tag{28}$$

In other words, the elements of $C^*_{n,l}$ are the elements of the power set of $\mathcal{N}_n$ with cardinality $l$. It is clear that such a set enumerates all the possible ways to pick $l$ integers from $n$ total integers ranging from 0 to $(n-1)$ without replacement. The cardinality of this set is necessarily equal to $_nC_l$.

This definition allows one to construct the elements of the set by carefully listing all the subsets of $\mathcal{N}_n$ while eliminating the ones that do not have a cardinality of $l$. What it does not specify is the order in which one may list the elements of $C^*_{n,l}$. After all, by definition, sets are unordered.

In the realms of probability and computer science, it is often necessary to list the elements of the set $C^*_{n,l}$. The question of the proper listing order invariably comes up. A few algorithms for achieving that goal have been listed by Knuth in his encyclopedic volumes of computer algorithms [15]. Many of these algorithms may be mapped to the following form.

$$C^{**}_{n,l} = \bigcup_{k_1=0}^{n-l} \bigcup_{k_2=k_1+1}^{n-l+1} \cdots \bigcup_{k_{l-1}=k_{l-2}+1}^{n-2} \bigcup_{k_l=k_{l-1}+1}^{n-1} \{k_1, k_2, \cdots, k_{l-1}, k_l\} \tag{29}$$

This form will be key to the remaining analysis in this article. Before using it, it is necessary to show that $C^*_{n,l}$ is equivalent to $C^{**}_{n,l}$. Notice that each element of $C^{**}_{n,l}$ contains $l$ integers. According to the limits listed in (29), it is clear that all of these integers are within the range $[0, n-1]$. Furthermore, the sequence $<k_1, k_2, \ldots, k_l>$ is listed in ascending order, because the lower limit of $k_{i+1}$ is $k_i + 1$. This also means that $k_i \neq k_{i'}$ as long as $i \neq i'$. Because each elements of $C^{**}_{n,l}$ is a set of $l$ unique natural numbers less than $n$, one may conclude that all elements of $C^{**}_{n,l}$ are also elements of $C^{*}_{n,l}$.

It is also possible to show that equation (29) contains $_nC_l$ "terms". If each "term" is unique, then the cardinality of $C^{**}_{n,l}$ is also equal to $_nC_l$. First, one may notice that the equation collapses to the following form when $l = 1$.

$$C^{**}_{n,1} = \bigcup_{k_1=0}^{n-1} \{k_1\} \tag{30}$$

As a result, the number of distinct elements in $C^{**}_{n,1}$ is equal to $_nC_1$, in accordance with the conjecture.

Consider another specific case, where $n = l$. Equation (29) reduces to the following form.

$$C^{**}_{n,n} = \{0, 1, \cdots, n-2, n-1\} \tag{31}$$

The cardinality of $C^{**}_{n,n}$ is $_nC_n = 1$, which also fits the conjecture. This is also sufficient to show that $|C^{**}_{2,l}| = {}_2C_l$.

For $n \geq 3$, one may attempt an inductive proof. Suppose $|C^{**}_{n-1,l}| = {}_{(n-1)}C_l$. The corresponding form for equation (31) is:

$$C^{**}_{n-1,l} = \bigcup_{k_1=0}^{n-l-1} \bigcup_{k_2=k_1+1}^{n-l} \cdots \bigcup_{k_{l-1}=k_{l-2}+1}^{n-3} \bigcup_{k_l=k_{l-1}+1}^{n-2} \{k_1, k_2, \cdots, k_{l-1}, k_l\} \tag{32}$$

Comparing the limits of the indices $k_1, k_2, \ldots, k_l$ in (32) and (29), it is clear that $C^{**}_{n-1,l}$ is a subset of $C^{**}_{n,l}$. Let $X$ be the part of $C^{**}_{n,l}$ with $C^{**}_{n-1,l}$ removed. The set $X$ must be in the form $\{k_1, k_2, \ldots, k_{l-1}, n-1\}$, since the integer $(n-1)$ is not available in the elements of $C^{**}_{n-1,l}$ and so must end up in the elements of $X$. On the flip side, an element of the form $\{k_1, k_2, \ldots, k_{l-1}, k'\}$, where $k' < (n-1)$, must belong in $C^{**}_{n-1,l}$, as all the values in the sequence clearly fall within the limits of the indices in (29).

$$C^{**}_{n,l} = C^{**}_{n-1,l} \cup X \tag{33}$$

$$X = \bigcup_{k_1=0}^{n-l} \bigcup_{k_2=k_1+1}^{n-l+1} \cdots \bigcup_{k_{l-1}=k_{l-2}+1}^{n-2} \{k_1, k_2, \cdots, k_{l-1}, n-1\} \tag{34}$$

Comparing (34) with (32), one may conclude that $X$ is equivalent to the set formed when the integer $(n-1)$ is appended to each and every element of $C^{**}_{n-1, l-1}$. The cardinality of $X$ is thus the same as that of $C^{**}_{n-1, l-1}$. This insight allows one to write:

$$|C^{**}_{n,l}| = |C^{**}_{n-1,l}| + |C^{**}_{n-1,l-1}| \tag{35}$$

Recall that the cardinalities of $C^{**}_{n-1, l}$ and $C^{**}_{n-1, l-1}$ are given to be $_{n-1}C_l$ and $_{n-1}C_{l-1}$, respectively. Using Pascal's identity, equation (35) becomes:

$$|C^{**}_{n,l}| = \binom{n}{l} \tag{36}$$

This completes the proof for the assertion that the cardinality of $C^{**}_{n, l}$ is equal to $_nC_l$. We have shown that all elements of $C^{**}_{n, l}$ are elements of $C^{*}_{n, l}$, and the cardinality of $C^{**}_{n, l}$ is equal to the cardinality of $C^{*}_{n, l}$. It must follow that there exists a bijection between $C^{**}_{n, l}$ and $C^{*}_{n, l}$.

Equation (32), in itself, is not a full description of an algorithm for listing all elements of $C^{*}_{n, l}$. This is because the order of the elements is still not fully specified. For the purpose of this article, the following algorithm will be used to clarify the proper order to be used.

1. Set the upper limits $\langle k_{U, 1}, k_{U, 2}, \ldots, k_{U, l}\rangle = \langle n-l, n-l+1, \ldots, n-1\rangle$.
2. Initialize $\langle k_1, k_2, \ldots k_l\rangle = \langle 0, 1, \ldots, l-1\rangle$.
3. Initialize $i = l$.
4. If $k_i \geq k_{U, i}$, decrement $i$ and repeat step 4 if the result is greater than or equal to 1. The algorithm ends if $i$ is now equal to 0.
5. Increment $k_i$.
6. Initialize $j = i + 1$.

7. Set $k_j = k_i + j - i$.
8. Increment $j$ by 1. If $j > l$, Output $<k_1, k_2, \ldots k_l>$ as the result for the current iteration and start the next iteration by going back to step 3. Return to step 7 otherwise.

This algorithm lists elements of $C^*_{n, l}$ in lexicographical order, with integers in each element listed in increasing order. For a concrete example, suppose $n = 5$ and $k = 3$. This algorithm results in the following ordering of the elements of $C^{**}_{5, 3}$. This list should be read from left to right, then top to bottom.

⟨0,1,2⟩  ⟨0,1,3⟩  ⟨0,1,4⟩  ⟨0,2,3⟩  ⟨0,2,4⟩
⟨0,3,4⟩  ⟨1,2,3⟩  ⟨1,2,4⟩  ⟨1,3,4⟩  ⟨2,3,4⟩

Because each iteration of the algorithm maps directly to an element in equation (29), each iteration produces a distinct element. One may show this independently by noticing that the steps in the algorithm described above are reminiscent of how one would perform addition to a base-$n$ integer. In fact, since all values in the sequence $<k_1, k_2, \ldots k_l>$ are within the range $[0, n - 1]$, one may define an injective mapping between such sequences and natural numbers:

$$F(\langle k_1, \ldots, k_l \rangle) = \sum_{i=0}^{l-1} k_{l-i} n^i \qquad (37)$$

As is the case for any pair of base-$n$ integers, the integer $F(<k_1, k_2, \ldots, k_l>)$ is larger than $F(<k'_1, k'_2, \ldots, k'_l>)$ if the first non-zero value in the sequence $<k_1 - k'_1, k_2 - k'_2, \ldots, k_l - k'_l>$ is larger than zero.

It is well-known that the "greater than" relationship for natural numbers is transitive, such that if $a > b$ and $b > c$, that $a > c$. Suppose $<k_1, k_2, \ldots k_l>$ is the output of a particular iteration $k$ of the algorithm, and $<k'_1, k'_2, \ldots k'_l>$ is the output of the following iteration ($k + 1$), one may map these two outputs to their corresponding integers $Z_k = F(<k_1, k_2, \ldots k_l>)$ and $Z_{k+1} = F(<k'_1, k'_2, \ldots k'_l>)$. If $Z_{k+1}$ is always greater than $Z_k$, one may leverage the transitivity of the "greater than" relationship of natural numbers to conclude that the outputs of the algorithm always "increase" monotonically. The immediate consequence is that no two different iterations produce the same output.

It is possible to show that $Z_{k+1}$ above is indeed larger than $Z_k$. If the output of a particular iteration $k$ of the algorithm is $<k_1, k_2, \ldots k_l>$, the output of the following iteration $(k + 1)$ is in a form $<k'_1, k'_2, \ldots k'_l>$ such that $k_1 = k'_1$, $k_2 = k'_2$, ..., $k_{i-1} = k'_{i-1}$ $(0 \leq i < l)$ remain unchanged, $k'_i = k_i + 1$, and $k'_{i+1}, k'_{i+2}, \ldots, k'_l$ are reset to the lowest valid values. If these two sequences are mapped to integers using equation (37), it is clear that the integer corresponding to the output of the iteration $k$ is smaller than that of the iteration $(k + 1)$ because the magnitudes of these two integers are determined by $k_i$ and $k'_i$, and $k'_i$ is always bigger. This completes the proof that no two iterations of the algorithm produce the same output.

From here on, the output of the $k^{th}$ iteration of the algorithm will be annotated as $C_{n, l, k}$. The sequence $<C_{n, l, 1}, C_{n, l, 2}, \ldots, C_{n, l, N}>$, where $N = {}_nC_l$, is a complete listing of all elements of $C^*_{n, l}$.

## Calculating the Higher Derivatives of the Falling Factorial

Finding the $l^{th}$ derivative of the falling factorial for a fixed, small value of $l$ is, no doubt, an elementary exercise of applying the Leibniz rule of differentiation repeatedly. The challenge is in coming up with a general solution for a variable $l$, and to annotate the result in a reasonably compact manner. We now have the tools to do both.

It is easy to verify that the following results for the second- and third derivatives of the falling factorial are correct.

$$F_n^{(2)}(x) = \begin{cases} 0 & 0 \leq n \leq 1 \\ 2 \sum_{k=1}^{\binom{n}{2}} \Theta(n, C_{n,2,k}; x) & n \geq 2 \end{cases} \quad x \in \mathbb{C} \quad (38)$$

$$F_n^{(3)}(x) = \begin{cases} 0 & 0 \leq n \leq 2 \\ 6 \sum_{k=1}^{\binom{n}{3}} \Theta(n, C_{n,3,k}; x) & n \geq 3 \end{cases} \quad x \in \mathbb{C} \quad (39)$$

It seems reasonable to make the following conjecture:

$$F_n^{(l)}(x) = \begin{cases} 0 & 0 \leq n \leq l-1 \\ \sum_{k=1}^{\binom{n}{l}} w_{n,l,k} \Theta(n, \mathcal{C}_{n,l,k}; x) & n \geq l \end{cases} \quad x \in \mathbb{C} \quad (40)$$

This equation is an assertion that the $l^{\text{th}}$ derivative of the falling factorial function is a weighted sum of the complete set of falling factorials of order $n$ with $l$ missing factors. There are $_nC_l$ ways to remove $l$ factors from the falling factorial of order $n$, and it is clear that equation (40) goes through all of these possibilities. The "weights", or coefficients, in the weighted sum are annotated as $w_{n, l, k}$.

One may prove this using induction. First, by comparing with equation (27), one can verify that (40) is true for $l = 1$. Suppose the equation is true for a particular value of $l$, then one may differentiate both sides to give:

$$F_n^{(l+1)}(x) = \sum_{k=1}^{\binom{n}{l}} w_{n,l,k} \frac{d}{dx}\{\Theta(n, \mathcal{C}_{n,l,k}; x)\} \quad (41)$$

According to equation (8), the derivative of the falling factorial with $l$ missing factors is a sum of $(n - l)$ different falling factorials with $(l + 1)$ missing factors. Each term of the result must retain the $l$ exact missing factors before the differentiation. One may then write:

$$F_n^{(l+1)}(x) = \sum_{k=1}^{\binom{n}{l}} w_{n,l,k} \sum_{k' \in (\mathcal{N}_n - \mathcal{C}_{n,l,k})} \Theta(n, \mathcal{C}_{n,l,k} \cup k'; x) \quad (42)$$

It is now possible to claim that the form of equation (42) allows one to state with confidence that it is a weighted sum of all members of the complete set of falling factorials with $(l + 1)$ missing factors. It is a basic tenet of combinatorics that the single operation of picking $(l + 1)$ values from $n$ total elements is the same as the composite operation of picking $l$ elements from the $n$ elements, followed by picking one element from the remaining $(n - l)$ elements. The outer summation of equation (42) represents all possible combinations of picking $l$ elements, while the inner summation represents all possible ways to pick the remaining one element. Thus the right hand side of (42) is a weighted sum of the complete set of falling factorials with $(l + 1)$ missing factors. In other words:

$$F_n^{(l+1)}(x) = \begin{cases} 0 & 0 \leq n \leq l \\ \sum_{k=1}^{\binom{n}{l+1}} w_{n,l+1,k} \Theta(n, C_{n,l+1,k}; x) & n \geq l + 1 \end{cases} \quad (43)$$

This proves that the relationship (40) holds for a particular value of $l$ as long as it is true for $(l - 1)$. Since it is known to be true for $l = 1$, it is true for all $l \geq 1$.

Up to now, the actual values of the coefficients $w_{n, l, k}$ have not been calculated. Numerically, there is strong indication that $w_{n, l, k}$ is a constant value independent of the value of $k$. This may be proven rigorously in the following manner.

Suppose the following is a true statement for a particular value of $l$:

$$w_{n,l,k} = l! \quad (44)$$

Equations (42) and (43) can then be rewritten as:

$$l! \sum_{k=1}^{\binom{n}{l}} \sum_{k' \in (\mathcal{N}_n - C_{n,l,k})} \Theta(n, C_{n,l,k} \cup k'; x) = \sum_{j=1}^{\binom{n}{l+1}} w_{n,l+1,j} \Theta(n, C_{n,l+1,j}; x) \quad (45)$$

One may solve for $w_{n, l + 1, j}$ by counting the number of ways $C_{n, l, k} \cup k'$ on the left hand side can become equivalent to $C_{n, l + 1, j}$ on the right hand side. This means $C_{n, l, k}$ is a subset of $C_{n, l + 1, j}$, and $k'$ is also a member of $C_{n, l + 1, j}$. This is equivalent to choosing $l$ objects from a total of $(l + 1)$ objects, and so there are $_{l + 1}C_l = l + 1$ ways this can happen. As a result, $w_{n, l + 1, j}$ is equal to $l!(l + 1) = (l + 1)!$. This allows one to write the solution for the $l^{\text{th}}$ derivative of the falling factorial function as:

$$F_n^{(l)}(x) = \begin{cases} 0 & 0 \le n \le l-1 \\ l! \sum_{k=1}^{\binom{n}{l}} \Theta(n, \mathcal{C}_{n,l,k}; x) & n \ge l \end{cases} \quad x \in \mathbb{C} \qquad (46)$$

## Solution of the $l^{th}$ Derivative of the Falling Factorial in "Harmonic" Form

If $x$ is not an integer in the range $[0, n-1]$, it is clear that the falling factorial power with missing factors may be calculated by taking the full falling factorial and dividing out the undesired factors.

$$\Theta(n, \mathcal{C}_{n,l,k}; x) = x^{\underline{n}} \prod_{j \in \mathcal{C}_{n,l,k}} \frac{1}{x-j} \quad x \in \mathbb{C} \wedge x \notin \mathcal{C}_{n,l,k} \qquad (47)$$

With this in mind, one may rewrite equation (46) as:

$$F_n^{(l)}(x) = \begin{cases} 0 & 0 \le n < l \\ l! \, x^{\underline{n}} \sum_{k=1}^{\binom{n}{l}} \prod_{j \in \mathcal{C}_{n,l,k}} \frac{1}{x-j} & n \ge l \end{cases} \quad \begin{array}{l} x \in \mathbb{C} \\ x \notin \mathcal{C}_{n,l,k} \end{array} \qquad (48)$$

It is now time to apply (29), which transforms equation (48) into the following form.

$$F_n^{(l)}(x) = \begin{cases} 0 & 0 \le n < l \\ l! \, x^{\underline{n}} \sum_{k_1=0}^{n-l} \sum_{k_2=k_1+1}^{n-l+1} \cdots \sum_{k_l=k_{l-1}+1}^{n-1} \prod_{j=1}^{l} \frac{1}{x-k_j} & n \ge l \end{cases} \qquad (49)$$

This may be seen as a "harmonic" form of equation (46).

### Introducing the Elementary Symmetrical Harmonic Sum

For the purpose of this article, the following quantity will be referred to as the "elementary symmetric harmonic sum". It can be used to simplify solutions of the form shown in equation (49). The various identities that can be derived using the definition will also be very helpful for revealing deeper insights of all quantities involved.

$$H_{n,l,r} = \begin{cases} 0 & l < 0 \\ 1 & l = 0 \\ \sum_{k=r+1}^{n-l+1} \dfrac{H_{n,l-1,k}}{k} & l \geq 1 \end{cases} \quad n - l \geq r \tag{50}$$

Equation (50) is true for all $n - l \geq r$. Otherwise the following convention should be followed.

$$H_{n,l,r} = 0 \quad n - l < r \tag{51}$$

For now, the parameters $n$, $l$, $r$ are all integers. In addition, $n$ and $r$ are natural numbers including zero.

The recursive definition (50) can be "unrolled" to produce a more intuitive representation of the quantity.

$$H_{n,l,r} = \sum_{k_1=r+1}^{n-l+1} \sum_{k_2=k_1+1}^{n-l+2} \cdots \sum_{k_l=k_{l-1}+1}^{n} \prod_{j=1}^{l} \frac{1}{k_j} = \sum_{r<k_1<k_2<\cdots<k_l\leq n} \frac{1}{k_1 k_2 \cdots k_l} \tag{52}$$

Thus the elementary symmetric harmonic sum is equivalent to the elementary symmetric polynomial where the $(n - r)$ variables are evaluated at the multiplicative inverses of $r + 1, r + 2, \ldots, n$.

$$H_{n,l,r} = e_l\left(\frac{1}{r+1}, \frac{1}{r+2}, \ldots, \frac{1}{n}\right) \quad n - r \geq l \tag{53}$$

It is easy to show from the definition that the elementary symmetric harmonic sum is equivalent to the harmonic number with $l = 1$ and $r = 0$. As a result, one may see the elementary symmetric harmonic sum as a generalization of the harmonic number.

$$H_{n,1,0} = \sum_{k=1}^{n} \frac{1}{k} = H_n \qquad (54)$$

Another special case is:

$$H_{n,l,n-l} = \frac{1}{n^{\underline{l}}} \quad 0 \le l \le n \qquad (55)$$

Working with the elementary symmetric harmonic sum is relatively easy, as it follows a few simple recursive relationships.

$$H_{n,l,r} = H_{n,l,r+1} + \frac{H_{n,l-1,r+1}}{r+1} \qquad (56)$$

$$H_{n,l,r} = H_{n-1,l,r} + \frac{H_{n-1,l-1,r}}{n} \qquad (57)$$

$$H_{n+1,l,r+1} - H_{n,l,r} = \left(\frac{1}{n+1} - \frac{1}{r+1}\right) H_{n,l-1,r+1} \qquad (58)$$

The first relationship can be proven through direct evaluation of the definition (50). The third relationship is simply a combination of the first and second relationships. The second relationship can be proven using equation (52).

$$\begin{aligned}
H_{n,l,r} &- H_{n-1,l,r} \\
&= \sum_{k_1=r+1}^{n-l+1} \sum_{k_2=k_1+1}^{n-l+2} \cdots \sum_{k_l=k_{l-1}+1}^{n} \prod_{j=1}^{l} \frac{1}{k_j} \\
&- \sum_{k_1=r+1}^{n-l} \sum_{k_2=k_1+1}^{n-l+1} \cdots \sum_{k_l=k_{l-1}+1}^{n-1} \prod_{j=1}^{l} \frac{1}{k_j}
\end{aligned} \qquad (59)$$

Notice that the first elementary symmetric harmonic sum has $_nC_l$ terms, while the second one has $_{(n-1)}C_l$ terms. The first elementary symmetric harmonic sum clearly contains all the terms found in the second elementary symmetric harmonic sum. As a result, the difference is the sum of all terms in the first elementary symmetric harmonic sum with $k_l = n$. This insight allows one to write:

$$H_{n,l,r} - H_{n-1,l,r} = \frac{1}{n} \sum_{k_1=r+1}^{n-l+1} \sum_{k_2=k_1+1}^{n-l+2} \cdots \sum_{k_{l-1}=k_{l-2}+1}^{n-1} \prod_{j=1}^{l-1} \frac{1}{k_j} \tag{60}$$

This then reduces to equation (57).

Equation (57) may be used to construct a matrix that can be used to look up values of the elementary symmetric harmonic sum. The steps for this construction are:

1. Define the row index as $n$ and the column index as $l$. The value of $r$ is fixed.
2. Populate the entries for $l = 0$, which represent the first column of the matrix, using the value 1.
3. Populate the diagonal $n = l$, starting from the upper left, using equation (55).
4. Populate the rest of the values between the first column and the diagonal using the recursive relationship (57).
5. All other values in the unfilled upper-right triangle are set to zero.

The recursive relationship (56) may also be used to construct a similar matrix with a fixed value of $n$.

The most common way to generalize the harmonic number is to raise each term to a particular power.

$$H_n^{(v)} = 1 + \frac{1}{2^v} + \frac{1}{3^v} + \cdots + \frac{1}{n^v} \tag{61}$$

This quantity is widely studied due to, in no small part, the connection to Riemann's zeta function. It is possible to generalize the elementary symmetric harmonic number in a similar way. The definition, special values, and recursive relationships of this quantity are shown below. These identities are very similar to those of the elementary symmetric harmonic sum, and may be derived using essentially identical techniques.

$$H_{n,l,r}^{(v)} = \begin{cases} 0 & l < 0 \\ 1 & l = 0 \\ \sum_{k=r+1}^{n-l+1} \frac{H^{(v)}(n, l-1, k)}{k^v} & l \geq 1 \end{cases} \quad n - l \geq r \tag{62}$$

$$H_{n,1,0}^{(v)} = \sum_{k=1}^{n} \frac{1}{k^v} = H_n^{(v)} \tag{63}$$

$$H_{n,l,n-l}^{(v)} = \frac{1}{(n_-^l)^v} \tag{64}$$

$$H_{n,l,r}^{(v)} = H_{n,l,r+1}^{(v)} + \frac{H_{n,l-1,r+1}^{(v)}}{(r+1)^v} \tag{65}$$

$$H_{n,l,r}^{(v)} = H_{n-1,l,r}^{(v)} + \frac{H_{n-1,l-1,r}^{(v)}}{n^v} \tag{66}$$

### Specific Solutions at Integral Values of $x \leq -1$ or $x \geq n$

If $x$ is an integer in the range $(-\infty, -1]$ or $[n, \infty)$, equation (49) does not suffer from the divided-by-zero problem, and can be simplified directly. One important fact to keep in mind about derivatives of the falling factorial function is that it is always either a symmetric or anti-symmetric function around the axis $x = (n - 1) / 2$. The following equations summarize this property.

$$x^{\underline{n}} = (-1)^n (n - 1 - x)^{\underline{n}} \tag{67}$$

From this, it is immediately obvious that the following is true.

$$F_n^{(l)}(m) = (-1)^{n-l} F_n^{(l)}(n - m - 1) \tag{68}$$

It is thus not necessary to solve for the case $x \leq -1$ if a solution for $x \geq n$ is available. Evaluating equation (49) at integral values of $x \geq n$ gives:

$$F_n^{(l)}(m) = l!\, m^{\underline{n}} \sum_{k_1=0}^{n-l} \sum_{k_2=k_1+1}^{n-l+1} \cdots \sum_{k_l=k_{l-1}+1}^{n-1} \prod_{j=1}^{l} \frac{1}{m - k_j} \quad \begin{array}{l} n \geq l \\ m \in \mathbb{Z} \\ m \geq n \end{array} \tag{69}$$

The elementary symmetric harmonic sum notation developed above allows this to be simplified to:

$$F_n^{(l)}(m) = l!\, m^{\underline{n}} H_{m,l,m-n} \quad \begin{array}{l} m \in \mathbb{Z} \\ m \geq n \end{array} \tag{70}$$

The solution for $m \leq -1$ follows trivially if equation (68) is applied, and will not be listed here.

## Specific Solutions at Integral Values of $0 \leq x < n$

Due to the divide-by-zero issue discussed above, equation (49) may not be used directly for numerical calculation for integral values of $x$ in the range $[0, n - 1]$. The poles can be canceled out ahead of time of the numerical calculation to create a more friendly expression. However, this procedure is not extremely straightforward, and some theorizing is needed to perform this task successfully in the general case.

Now that $0 < x < n$ is taken as a given, it is clear that the factor $(x - k_j)$ in (49), where $0 \leq k_j < n$, is going to be a value $y$ within the range $[-(n - x), x]$. Suppose $(n - x) > x$, then for every value of $y$ where $1 \leq y \leq x$, there is a $-y$ within the range $[-(n - x), -1]$. If both $y$ and $-y$ are in found within the range $[-(n - x), x]$, one may say that the value $y$ is "paired". For those values of $-y$ where $-(n - x) \leq -y < -x$, the corresponding positive value $y$ is not in the range $[-(n - x), x]$. One may say that these values of $-y$ are "unpaired".

The same argument can be applied on the case where $(n - x) < x$. In this case, for every value of $-y$ where $-(n - x) \leq -y \leq -1$, there is a paired $y$ within the range $[1, x]$. The values of $y$ with the range $[(n - x + 1), x]$ are unpaired.

As a concrete example, suppose $n = 5$, such that the falling factorial is $x(x - 1)(x - 2)(x - 3)(x - 4)$. Evaluating this at $x = 1$ gives the product $1 \times 0 \times -1 \times -2 \times -3$. In this case, the factors $x$ and $(x - 2)$ are "paired" because their absolute values are the same, and the factors $(x - 3)$ and $(x - 4)$ are unpaired.

For the purpose of the current discussion, suppose the set of values with a paired negative counterpart is defined as set A, while the unpaired values are in the set B. Let C be the superset containing all values in sets A and B.

Recall that the value of *l* signifies the number of factors in the denominator within each term in (49). For each term, one picks *l* distinct values from set C. From a combinatorial point of view, this is equivalent to picking *j* values from set A and (*l* − *j*) values from set B.

Let *j* = 1. Because the value picked from set A has a negative counterpart, it is always possible to find another term that is the negative of this term. In other words, it is not necessary to evaluate terms where *j* = 1, as it will eventually be canceled out to zero.

The same situation applies for any odd values of *j*. There is always a term that is the negative of the current term.

Let *j* = 2. If the absolute values of the two elements picked from set A are different, then it is always possible to pick another two values with the exact same product but a different sign. However, if the absolute values of the two elements picked are the same, there will not be another term with an opposite sign to cancel it out. In fact, for any even value of *j*, there are un-canceled terms if the elements are picked in pairs of the same absolute values.

This argument leads to the following sanitized version of equation (49) for $0 \leq x < n$.

$$F_n^{(l)}(m) = (-1)^{n-m-1} l! \, m! \, (n - m - 1)! \times \sum_k (-1)^k H_{n-m-1,k,0}^{(2)} H_{m,l-2k-1,n-m-1} \qquad \begin{array}{c} m \in \mathbb{Z} \\ \frac{n-1}{2} \leq m < n \end{array} \qquad (71)$$

This completes the derivation of the solutions for the higher derivatives of the falling factorial evaluated at integral values.

### Relationship between the Elementary Symmetric Harmonic Sum and the Stirling Numbers of the First Kind

It was noted before in equation (53) that the elementary symmetric harmonic sum is, in fact, an elementary symmetric polynomial evaluated at specific values. Also, |*s*(*n*, *l*)| is itself equal to $e_{n-l}$ (1, 2, ..., *n* − 1). In other words, |*s*(*n* + 1, *l* + 1)| is

equal to $e_{n-l}(1, 2, \ldots, n)$. Dividing this by $n!$ gives the following simple relationship between the "complete" elementary symmetric harmonic sum and Stirling numbers of the first kind:

$$H_{n,l,0} = \frac{|s(n+1, l+1)|}{n!} \quad 0 \leq l \leq n \tag{72}$$

If one studies the elementary symmetric harmonic sum numerically, it will appear that any elementary symmetric harmonic sum may be written as a weighted sum of Stirling numbers of the first kind. This proposition may be written as:

$$H_{n,l,r} = \frac{(-1)^{n+l}}{n!} \sum_{k=0}^{l+1} A_{r,l-k+1} s(n+1, k) \tag{73}$$

$$A_{r,k} = \begin{cases} 1 & k = 0 \\ 0 & r = 0, k \geq 1 \\ A_{r-1,k} + \frac{A_{r,k-1}}{r} & r \geq 1, k \geq 1 \end{cases} \tag{74}$$

This relationship may also be proven through induction. The initial conditions are $r = 0$ ($k$ arbitrary) and $k = 0$ ($r$ arbitrary). The first initial condition has already been proven, as shown in equation (72). The second condition is taken as a given and will be shown to be necessary shortly.

First, one may show that the proposition in equation (74) follows the correct recursive relationship. Plugging it into the proven relationship (56) gives:

$$\sum_{k=1}^{l+1} A_{r,l-k+1} s(n+1, k) = \sum_{k=1}^{l+1} \left( A_{r+1,l-k+1} - \frac{A_{r+1,l-k}}{r+1} \right) s(n+1, k) \tag{75}$$

This equality is true if the definition (74) is used. Notice that $A_{r,k}$ must be equal to 1 for the equality to hold, as required. This completes the proof for (73) and (74).

Equation (73) has the following special case when $r = 1$.

$$\sum_{k=1}^{l+1} s(n+1,k) = (-1)^{n+l} n! \, H_{n,l,1} \tag{76}$$

This gives a simple expression for the incomplete sum of a row of the Stirling numbers of the first kind. This would have been relatively difficult to describe without the elementary symmetric harmonic sum notation.

The values for $A_{r,k}$ may be listed as a matrix with $r$ as the row index and $k$ as the column index. This matrix as no limit is either the row or column direction. The values in the first column is always 1, and the values in the second column are the harmonic numbers. The value of $A_{r,k}$ appears to tend to $r$ as $k$ tends to infinity. This will be apparent later.

If one attempts to write out the first few values of $A_{r,k}$, one may notice a tantalizing pattern. It appears that the following is a true identity.

$$A_{r,k} = \sum_{j=1}^{r} \frac{(-1)^{j+1}}{j^k} \binom{r}{j} \tag{77}$$

It can be proven to be true simply by plugging this conjecture back into equation (74).

This identity is interesting because it is very close to the definition of the Stirling numbers of the second kind. The major difference is that the "row" parameter is a negative number. This type of extended Stirling number of the second kind was studied by Branson [16], who called it the negative-positive Stirling number of the second kind. This realization allows one to rewrite equation (73) as:

$$H_{n,l,r} = \frac{(-1)^{n+1}}{n^{n-r}} \sum_{k=0}^{l+1} S(-l+k-1,r) s(n+1,k) \tag{78}$$

$$A_{r,k} = (-1)^{r+1} r! \, S(-k,r) \tag{79}$$

This is somewhat of a surprising result. It is interesting to note that Branson discovered that the absolute value of the negative-positive Stirling number of the second kind is the same as that of the negative-positive Stirling number of the first kind when the column and row indices are switched. Thus (78) can be written in terms of negative-positive Stirling numbers of the first kind as well.

$$H_{n,l,r} = \frac{(-1)^{n+1}}{n^{\underline{n-r}}} \sum_{k=0}^{l+1} (-1)^k s(-r, l-k+1) s(n+1, k) \tag{80}$$

Before going further, it is important to examine some of the properties of these negative-positive Stirling numbers. In [1], Knuth provided an overview of what Branson called negative-negative Stirling numbers. It was noted that the Stirling numbers of the first and second kinds, are, in fact, two representations of the same underlying quantity. This "law of duality" is summarized with the following identity.

$$S(n, k) = (-1)^{n+k} s(-k, -n) \tag{81}$$

Knuth noted that the classic recursive relationships for Stirling numbers are valid for all positive, negative, and zero values of $n$ and $k$ if the negative-positive Stirling numbers are set to zero through proper definition of the boundary conditions. The problem with this approach of extended the Stirling numbers to negative parameters is that the connection to the following well-known formula is lost.

$$S(n, k) = \frac{1}{k!} \sum_{j=1}^{k} (-1)^{k-j} \binom{k}{j} j^n + [n = k = 0] \tag{82}$$

It is quite clear that this definition may be extended to negative values of $n$ and positive values of $k$ with non-zero results. This is exactly the premise of equation (77). However, this is no way to preserve equation (82) and the recursive relationship for Stirling numbers at the same time. In fact, to generate the negative-positive Stirling numbers, the following boundary condition must be used.

$$s(n, 0^+) = \frac{[n \leq 0]}{(-n)!} \tag{83}$$

Whereas, to generate negative-negative Stirling numbers, the following boundary condition must be used.

$$s(n, 0^-) = [n = 0] \tag{84}$$

Therefore, the proper value of $s(n, 0)$ to use depends on the "direction" of the summation. In all cases below, we will be using the positive k direction, and so the definition (83) will be used. To prevent any confusion, a table of the extended values of $s(n, k)$ are shown below.

| n / k | 0 | 1 | 2 | 3 | 4 | 5 |
|---|---|---|---|---|---|---|
| -5 | 1.39e-3 | -1.90e-2 | 2.78e-2 | -3.38e-2 | 3.73e-2 | -3.94e-2 |
| -4 | 8.33e-3 | -8.68e-2 | 0.120 | -0.141 | 0.153 | -0.160 |
| -3 | 4.17e-2 | -0.306 | 0.394 | -0.444 | 0.471 | -0.485 |
| -2 | 0.167 | -0.75 | 0.875 | -0.938 | 0.969 | -0.984 |
| -1 | 0.5 | -1 | 1 | -1 | 1 | -1 |
| 0 | 1 | 0 | 0 | 0 | 0 | 0 |
| 1 | 0 | 1 | 0 | 0 | 0 | 0 |
| 2 | 0 | -1 | 1 | 0 | 0 | 0 |
| 3 | 0 | 2 | -3 | 1 | 0 | 0 |
| 4 | 0 | -6 | 11 | -6 | 1 | 0 |
| 5 | 0 | 24 | -50 | 36 | -10 | 1 |

If one defines the positive-positive Stirling numbers in such a way that they are always sign-less, the negative-positive version of those numbers can still be negative. Furthermore, the sign difference between positive-positive $s(n, r)$ and negative-negative $S(-n, -r)$ is $(-1)^{n+r}$, while the sign difference between negative-positive $s(-n, r)$ and $S(-r, n)$ is $(-1)^{n+r+1}$. The author found that using the sign-less version of the Stirling numbers can actually cause more confusion.

With this understanding of the elementary symmetric harmonic sums and extended Stirling numbers, it is now possible to write the solutions of the higher derivatives of the falling factorial in a new form.

$$F_n^{(l)}(m) = l! \sum_{k=0}^{l+1} (-1)^{m+k+1} s(n-m, l-k+1) s(m+1, k) \quad \begin{array}{c} m \in \mathbb{Z} \\ m \geq n \end{array} \tag{85}$$

$$F_n^{(l)}(m) = (-1)^n l! \sum_{j=0}^{l}(-1)^j s(m+1,j)$$
$$\sum_{k=0}^{\lfloor\frac{l-j}{2}\rfloor}(-1)^k s(-n+m+1, l-2k-j)H^{(2)}_{n-m-1,k,0}$$
$$m \in \mathbb{Z}$$
$$\frac{n-1}{2} \leq m < n \quad (86)$$

The second equation can be simplified vastly using the following identity:

$$(-1)^{l+j+1}(m!)^2 \sum_{k=0}^{m}(-1)^k s(-m, l-j-2k)H^{(2)}_{m,k,0} = s(m+1, l-j+1) \quad (87)$$

This identity can be proven by showing that the left hand side follows the same recursive relationship as the Stirling number of the first kind. This can be done by applying equation (74) on the left hand side twice, followed by reducing the result using the recursive relationship (66) of the generalized elementary symmetric harmonic sum. The full work will not be shown here.

The identity allows one to simplify (86) to the following.

$$F_n^{(l)}(m) = (-1)^{m+1} l! \sum_{j}(-1)^j s(n-m, l-j+1)s(m+1, j) \quad \begin{array}{c} m \in \mathbb{Z} \\ \frac{n-1}{2} \leq m < n \end{array} \quad (88)$$

It is interesting to note that (88) has the same form has (85), which is the corresponding solution for $m \geq n$. It is a simple exercise using equation (67) to show that this solution, in fact, works for all integral values of $m$. Therefore, one may write the following single-line solution for the higher derivatives of the falling factorial.

$$F_n^{(l)}(m) = l! \sum_{k=0}^{l+1}(-1)^{m+k+1} s(n-m, l-k+1)s(m+1, k) \quad m \in \mathbb{Z} \quad (89)$$

An equivalent form, which, again, is a consequence of equation (67), is shown below.

$$F_n^{(l)}(m) = l! \sum_{k=0}^{l+1} (-1)^{m+l+k} s(m+1, l-k+1) s(n-m, k) \quad m \in \mathbb{Z} \quad (90)$$

This is made possible only by allowing the use of negative-positive Stirling numbers.

## Some Interesting Identities Based on Higher Derivatives of the Falling Factorial

Consider the Taylor expansion of the falling factorial power around $x = m$:

$$x^{\underline{n}} = \sum_{l=0}^{n} \frac{(x-m)^l}{l!} F_n^{(l)}(m) \quad (91)$$

Taking the derivative of both sides of this equation gives the following result.

$$\frac{1}{l'!} F_n^{(l')}(m') = \sum_{l=l'}^{n} \binom{l}{l'} \frac{1}{l!} F_n^{(l)}(m)(m'-m)^{l-l'} \quad (92)$$

At this point, it is useful to summarize the results obtained in this article thus far. We now have solutions for the higher derivatives of the falling factorial for all integral values of $x$. The solutions takes different forms depending on which of the following regions $x$ falls within:

1. "Harmonic" and "Stirling" forms for $x \leq -1$.
2. "Harmonic" form for $0 \leq x \leq (n-1)/2$.
3. "Harmonic" form for $(n-1)/2 \leq x < n$.
4. "Stirling" form for $-1 \leq x \leq n$.
5. "Harmonic" and "Stirling" forms for $x \geq n$.

One may randomly pick any two from above to substitute into (92). Since there are 8 different equations per slot, this exercise can potentially create 64 different identities. Some of the more interesting ones will be presented below.

$$\sum_{k \geq 0}(-1)^k s(n-m+1,k) \sum_{l \geq k-1} \binom{l}{l'} s(m, l-k+1) r^{l-l'} =$$
$$(-1)^{r+l'+1} \sum_{k=0}^{l'+1}(-1)^k s(n+r-m+1, l'-k+1) s(m-r, k) \quad (93)$$

The right side hand can also be written in its "harmonic" form with equation (80). By evaluating this identity at certain specific values, one gets a large number of more specialized identities, most of which have been discussed in other works on Stirling numbers identities. A couple examples are shown below.

$$\sum_{k \geq 0} \binom{k}{l'} s(n+1, k+1) n^{k-l'} = |s(n, l')| \quad (94)$$

$$\sum_{k \geq 0} \binom{k}{l'} s(n+1, k+1) 2^{k-l'} = s(n-1, l') + s(n-1, l'-1) \quad (95)$$

A particular form of (93) may be used as the means for extending the elementary symmetric harmonic sum into real or complex values of $r$.

$$H_{n,l',0}(r) = \frac{1}{(n+r)_n} \sum_{k \geq 0} \binom{k}{l'} |s(n+1, k+1)| r^{k-l'} \quad r \in \mathbb{C} \quad (96)$$

$$H_{n,l',m}(r) = H_{n+m+r,l',m+r} \quad r \in \mathbb{N} \quad (97)$$

Another interesting set of identities that can be generated from equation (92) is shown below.

$$\sum_{l=l'}^{n} \binom{l}{l'} H_{n+m,l,m} r^{l-l'} =$$
$$\frac{(-1)^{n+m+l'+r}}{(n+m)^{\underline{n}}} \sum_{k=0}^{l'+1}(-1)^k s(n+m+r+1, l'-k+1) s(-m-r, k) \quad (98)$$

Again, the right hand side can be changed to its "harmonic" form with equation (80). This gives us a similar equation to (96).

$$H_{n,l',m}(r) = \frac{(n+m)^{\underline{n}}}{(n+m+r)^{\underline{n}}} \sum_{l=l'}^{n} \binom{l}{l'} H_{n+m,l,m} r^{l-l'} \quad r \in \mathbb{C} \tag{99}$$

Evaluating equation (98) at specific values gives the following identities.

$$\sum_{l=l'}^{n} \binom{l}{l'} H_{n+m,l,m} = \frac{n+m+1}{m+1} H_{n+m+1,l',m+1} \tag{100}$$

$$\sum_{l=0}^{n} H_{n+m,l,m} = \frac{n+m+1}{m+1} \tag{101}$$

$$\sum_{l=0}^{n} (-1)^l H_{n+m,l,m} = \begin{cases} \frac{m}{n+m} & m \geq 1 \\ \delta_n & m = 0 \end{cases} \tag{102}$$

$$\sum_{l=l'}^{n} \binom{l}{l'} H_{n+m,l,m} (-m)^{l-l'} = \frac{|s(n+1, l'+1)|}{(n+m)_n} \tag{103}$$

$$\sum_{l=0}^{n} H_{n+m,l,m} (-m)^l = \binom{n+m}{n}^{-1} \tag{104}$$

$$\sum_{l=0}^{n} H_{n+m,l,m} (-m-1)^l = \delta_n \tag{105}$$

### Tying Up Loose Ends – Symmetric Functions and *r*-Stirling Numbers

As mentioned earlier, the Stirling numbers are strongly related to symmetric polynomials. In particular, $|s(n + 1, n + k - 1)| = e_k (1, 2, …, n)$ for positive-positive Stirling numbers of the first kind. The *r*-Stirling numbers of the first kind [12] may be seen as a generalization of $|s(n, k)|$ through this relationship:

$$\left[ \begin{array}{c} n+1 \\ n-k+1 \end{array} \right]_r = e_k(r, r+1, \cdots, n) \quad n-r+1 \geq k \tag{106}$$

In comparison, the elementary symmetric harmonic sum is the same as $H_{n,k,r-1} = e_k(1/r, 1/(r+1), \ldots, 1/n)$. As a result, the $r$-Stirling numbers are connected to the elementary symmetric harmonic number through the following equation:

$$H_{n,l,r} = \frac{1}{n^{\underline{n-r}}} \begin{bmatrix} n+1 \\ l+r+1 \end{bmatrix}_{r+1} \qquad n - l \geq k \tag{107}$$

Using equation (70), one gets:

$$F_n^{(l)}(m) = l! \begin{bmatrix} m+1 \\ l+m-n+1 \end{bmatrix}_{m-n+1} \tag{108}$$

This is consistent with Broder's results.

Due to the form of equation (111) above, any summation involving the elementary symmetric harmonic number with $l$ as the index can be directly translated into a valid identity involving $r$-Stirling numbers. As a result, all of the identities (99-105) have corresponding forms using $r$-Stirling numbers.

The "$r$-Stirling polynomials" described in [12] may be written as:

$$H_{n,l',m}(r) = \frac{1}{(n+m+r)^{\underline{n}}} R_1(n, l', m+r+1) \tag{109}$$

From the symmetric sum definition of the $r$-Stirling numbers, it is clear that these numbers are related to the coefficients of the falling factorials with missing factors described previously. In particular, the $r$-Stirling numbers are the coefficients of a falling factorial with the first $r$ factors removed.

Another specific type of symmetric polynomials is called the complete homogeneous symmetric polynomials. This is defined as:

$$h_k(X_1, X_2, \ldots, X_n) = \sum_{1 \leq j_1 \leq j_2 \leq \cdots \leq j_k \leq n} X_{j_1} X_{j_2} \cdots X_{j_n} \tag{110}$$

The positive-positive Stirling numbers of the second kind and the negative-positive Stirling numbers may be written in terms of this type of symmetric polynomial.

$$S(n+k, k) = h_n(1, 2, \ldots, k) \tag{111}$$

$$(-1)^{k-1} k!\, S(-n, k) = (-1)^n k!\, s(-k, n) = h_n\left(1, \frac{1}{2}, \ldots, \frac{1}{k}\right) \tag{112}$$

Therefore, many of the quantities encountered in this note so far may be seen as a coherent whole. When the harmonic numbers and the Stirling numbers of the first kind are generalized by manipulating the elementary symmetric polynomial, they form the elementary symmetric harmonic sum and the $r$-Stirling numbers, which are trivially related to each other through equation (107). When $|s(n, k)|$ is extended to negative $n$ and negative $k$ values, the result is a complete homogenous symmetric sum, which is equal to $S(k, n)$. The harmonic version of this complete homogeneous symmetric sum becomes the quantity $A_{r, k}$, described in equations (73-74).

## Conclusion

After applying a brute-force approach to calculating the higher derivatives of the falling factorial function, it was found a number of different quantities can be used to describe the results. In general, the solution is a weighted sum of all possible falling factorials with $l$ missing factors, where $l$ is the order of the derivative. When the falling factorial with missing factors is expanded in standard polynomial form, the coefficients are discovered to follow the same recursive relationship as the Stirling numbers of the first kind.

When the solution is evaluated at integral values, it was found that the result may be described using elementary symmetric harmonic sums. These sums are shown to have an interesting relationship with Stirling numbers, especially when the Stirling numbers are extended to allow negative parameters. They were also found to be related to $r$-Stirling numbers. A number of new identities involving these elementary symmetric harmonic sums and extended Stirling numbers were derived using the solutions of the higher derivatives of the falling factorial.

# References

[1] D. E. Knuth, "Two Notes on Notation," *The American Mathematical Monthly*, Volume 99, Number 5, May 1992, pp. 403-422.

[2] C. Jordan, Calculus of Finite Differences, Budapest, 1939.

[3] J. Stirling, *Methodus Differentialis*, 1730 (Transl.: in I. Tweddle, *Methodus Differentialis: An Annotated Translation of Stirling's Text*, London: Springer-Verlag, 2003).

[4] C. Tweedle, *James Stirling a Sketch of His Life and Works along with His Scientific Correspondence*, Oxford: Clarendon Press, 1922.

[5] N. Nielsen, *Handbuch der Theorie der Gammafunktion*, Leipzig: B. G. Teubner, 1906.

[6] J. Ginsburg, "Note on Stirling's Numbers," *The American Mathematical Monthly*, Vol. 35, No. 2, 1928, pp. 77-80.

[7] F. Viète, In Artem Analyticem Isagoge, Tours, 1591 (Transl.: in T. R. Witmer, *The Analytic Art*, Kent: Kent State Univ. Press, 1983).

[8] A. Girard, *Invention Nouvelle en l'Algebre*, 1629.

[9] L. Carlitz, "Weighted Stirling numbers of the first and second kind – I," *The Fibonacci Quarterly*, 18, 1980, pp. 147-162.

[10] L. Carlitz, "Weighted Stirling numbers of the first and second kind – II," *The Fibonacci Quarterly*, 18, 1980, pp. 242-257.

[11] M. Koutras, "Non-Central Stirling Numbers and Some Applications," *Discrete Mathematics*, vol. 42, 1982, pp. 73-89.

[12] A. Broder, "The r-Stirling Numbers," *Discrete Mathematics*, vol. 49, 1984, pp.241-259.

[13] I. Newton, *Mathematical Papers – Vol. III*, D. T. Whiteside ed., London: Cambridge Univ. Press, 1969.

[14] G. E. Andrews and K. Uchimura, "*Identities in Combinatorics IV: Differentiation and Harmonic Numbers*," Utilitas Mathematica, vol. 28, 1985, pp. 265-269.

[15] D. E. Knuth, *The Art of Computer Programming, Volume 4A: Combinatorial Algorithms, Part 1*, Addison-Wesley, 2011.

[16] D. Branson. "An extension of Stirling numbers," *The Fibonacci Quarterly*, Vol. 34, pp. 213-223.